\documentclass[11pt]{article}
\usepackage[margin=1in]{geometry}

\usepackage[pdftex]{graphicx}
\graphicspath{{../pdf/}{../jpeg/}}
\DeclareGraphicsExtensions{.pdf,.jpeg,.png}
\usepackage{amssymb,amsmath,url}
\usepackage{bm} 
\usepackage{multirow}
\usepackage{booktabs}	
\usepackage{epstopdf}
\usepackage{color}
\usepackage{subfigure}
\usepackage{bigfoot}		
\usepackage{algpseudocode}
\usepackage{pifont}
\usepackage{tikz}
\usepackage{color}
\usepackage{cite}
\usepackage{lipsum}
\usepackage{algorithm} 
\usepackage[pdftex]{graphicx}

\algrenewcommand\algorithmicrequire{\textbf{Input:}}
\algrenewcommand\algorithmicensure{\textbf{Output:}}
\algrenewcommand\algorithmicforall{\textbf{For}}

\newtheorem{theorem}{Theorem}
\newtheorem{lemma}{Lemma}
\newtheorem{proposition}{Proposition}

\allowdisplaybreaks
\newtheorem{definition}{Definition}

\newtheorem{remark}{Remark}

\newtheorem{assumption}{Assumption}
\DeclareMathOperator*{\argmin}{arg\,min}
\DeclareMathOperator*{\argmax}{arg\,max}

\title{\LARGE
	Safe Zeroth-Order Optimization Using Linear Programs
}
\author{Baiwei Guo\thanks{B. Guo and G. Ferrari-Trecate are with the DECODE group, Institute of Mechanical Engineering, EPFL, Switzerland. e-mails: {\tt \{baiwei.guo, giancarlo.ferraritrecate\}@epfl.ch}.} \and Yang Wang\thanks{Y. Wang is with Delft Center for Systems and Control, Delft
University of Technology. {\tt \{y.wang-40\}@tudelft.nl}.}\and Yuning Jiang \thanks{Y. Jiang is with the PREDICT group, Institute of Mechanical Engineering, EPFL, Switzerland. e-mails: {\tt yuning.jiang@epfl.ch}.} \and Maryam Kamgarpour\thanks{ M. Kamgarpour is with the SYCAMORE group, Institute of Mechanical Engineering, EPFL, Switzerland. e-mails: {\tt maryam.kamgarpour@epfl.ch}.} \and Giancarlo Ferrari-Trecate\footnotemark[1] }

\allowdisplaybreaks
\begin{document}
\maketitle \let\thefootnote\relax\footnotetext{This work was supported by the Swiss National Science Foundation under the NCCR Automation (grant agreement 51NF40\textunderscore 80545).}\\

\begin{abstract}\noindent To solve unmodeled optimization problems with hard constraints, this paper proposes a novel zeroth-order approach called Safe Zeroth-order Optimization using Linear Programs (SZO-LP). The SZO-LP method solves a linear program in each iteration to find a descent direction, followed by a step length determination. We prove that, under mild conditions, the iterates of SZO-LP have an accumulation point that is also the primal of a KKT pair. We then apply SZO-LP to solve an Optimal Power Flow (OPF) problem on the IEEE 30-bus system. The results demonstrate that SZO-LP requires less computation time and samples compared to state-of-the-art approaches.  
\end{abstract}
\section{Introduction}

A variety of applications, including power network operations \cite{chu2021frequency}, machine learning \cite{chen2020DFO}, trajectory optimization \cite{manchester2016derivative} and optimal control \cite{rao1998application}, require solving complex optimization problems with hard safety constraints. However, it is not always possible to obtain the expressions of the objective and constraint functions, or sufficient data on feasible system trajectories for modeling. In this context, safe zeroth-order optimization methods can be used to address unmodeled optimization problems with constraints. These methods rely solely on sampling by evaluating unknown objective and constraint functions at selected points \cite{Bajaj2021} and the term ``safe" refers to the feasibility of the samples (i.e., the satisfaction of the constraints).

Prominent safe zeroth-order methods include SafeOPT and its variations \cite{sui2015safe, turchetta2019safe, sabug2022smgo}. These approaches assume knowledge of a Lipschitz constant of the objective and constraint functions, while \cite{vinod2022constrained} utilizes a Lipschitz constant of function gradients (the smoothness constants). By using these quantities, one can build local proxies for the constraint functions. Starting from a feasible point, \cite{sui2015safe,sabug2022smgo,vinod2022constrained} utilize these proxies to search for potential minimizers. However, for each search, one has to use a global optimization method to solve a non-convex subproblem, which makes the algorithm computationally intractable for problems with many decision variables.

To reduce the computational complexity, another research direction involves incorporating barrier functions in the objective to penalize proximity to the boundary of the feasible set \cite{lewis2002globally,audet2009progressive}. The Extremum-Seeking methods \cite{hazeleger2022sampled} and the LB-SGD algorithm \cite{usmanova2022log} minimize a cost equipped with log-barrier penalty terms based on the estimated gradient. Although they do not require solving optimization subproblems, the performance of these two methods might not be satisfactory due to the log penalties. In Extremum Seeking, it can be challenging to tune the weight of the penalty term because a large weight can lead to suboptimality while a small weight might result in infeasibility. In LB-SGD, large values of the log barrier term and its derivative, when the iterates approach the boundary of the feasible set, can result in small step lengths and slow down convergence. 

Another approach to safe zeroth-order optimization is SZO-QQ proposed in \cite{guo2023safe}. It avoids log barrier penalties while still ensuring sample feasibility and is more sample-efficient than LB-SGD \cite{guo2023safe}. This is accomplished by utilizing convex quadratic proxies for the constraint functions to construct local feasible sets, over which the proxy for the objective function is then minimized. Unlike SafeOPT, the subproblems of SZO-QQ are convex Quadratically Constrained Quadratic Programs (QCQPs), which can be solved much faster than the non-convex subproblems in SafeOPT. However, SZO-QQ falls behind LB-SGD and Extremum Seeking in terms of computational efficiency (see Section~\ref{sec: numericals}) when dealing with large problems (with hundreds of constraints) because the size of each QCQP subproblem is almost the same as the original problem. In this paper, we propose a novel, safe zeroth-order method whose subproblems have much fewer constraints and can be computationally efficient.

Optimal Power Flow (OPF) is an example of large-scale optimization problems that can benefit from zeroth-order optimization. Its objective is to allocate the active and reactive power generation, transmission line flows and voltage levels to minimize costs while satisfying operational and security constraints such as transmission line capacity and voltage level limits. In recent years, OPF has gained considerable attention due to the rising demand for efficient and reliable operation of power systems, as well as the integration of renewable energy sources and energy storage systems \cite{abdi2017review}. However, the application of OPF to power system operation is a significant challenge due to the difficulties in accurately deriving a system model.  Therefore, we consider applying our model-free method to solve OPF problems.

The contributions of this paper are summarized as follows:
\begin{itemize}
    \item We present a novel  approach called Safe Zeroth-Order optimization using Linear Programs (SZO-LP). This method iteratively solves linear programming subproblems to derive descent directions and then decides the step length by sampling;
    \item We show that, under mild assumptions, a subsequence of SZO-LP's iterates converges to the primal of a KKT pair (see Definition \ref{def: KKT});
    \item By application to an IEEE 30-bus benchmark problem, we show that SZO-LP can efficiently solve an OPF problem with 11 decision variables and 158 constraints. We compare SZO-LP with state-of-the-art approaches and demonstrate its advantages in terms of computation time and the number of samples required.
\end{itemize}

\textit{Notations:}
We use $e_i\in\mathbb R^d$ to define the $i$-th standard basis of vector space $\mathbb{R}^d$ and $\|\cdot\|$ to denote the two norms throughout the paper. Given a vector $x\in \mathbb{R}^d$ and a scalar $\epsilon>0$, we write $x = [x^{(1)}, \ldots, x^{(d)}]^\top$ and $\mathcal B_\epsilon(x)=\{y:\|y-x\|\leq \epsilon\}$. We use $\mathbb{Z}_i^j=\{i,i+1,\ldots,j\}$ to define the set of integers ranging from $i$ to $j$ with $i<j$. For two vectors $x,y\in\mathbb{R}^d$, we use $\langle x,y\rangle :=x^\top y$ to define the inner product.

\section{Problem Formulation}
\label{sec: Prob}
We consider the constrained optimization problem
\begin{equation}
\label{eq: optimization problem}
\min_{x\in \mathbb{R}^d}  
\;\; f_0(x)\quad \text{subject to}\;\;
x\in \Omega,
\end{equation}
where $\Omega:=\{x: f_i(x)\leq 0, i\in\mathbb Z_1^m\}$ is the feasible set. The functions $f_i:\mathbb R^d\to\mathbb R$, $i\in\mathbb{Z}_0^m$, are unknown but can be sampled at query points. Throughout this paper, we make the following assumptions on the smoothness of the objective and constraint functions, availability of a strictly feasible point $x_0$ and boundedness of a sublevel set that includes~$x_0$.
\begin{assumption}
\label{ass: smoothness}
The functions $f_i(x)$, $i\in\mathbb Z_0^m$ are continuously differentiable and there are known constants $L_i,M_i>0$ such that for any $x_1$, $x_2\in \Omega $,
\begin{subequations}
\label{eq: smoothness}
\begin{align}
\label{eq: smoothness_1}
|f_i(x_1)-f_i(x_2)|&\leq L_i\|x_1-x_2\|,\\
\label{eq: smoothness_2}
\|\nabla f_i(x_1)-\nabla f_i(x_2)\|&\leq M_i\|x_1-x_2\|.
\end{align}
\end{subequations}
We assume $L_i>{\inf}\{L_i: \text{ \eqref{eq: smoothness_1}} \text{ holds},\forall x_1,x_2\in \Omega\}$ and $M_i>{\inf}\{M_i: \text{ \eqref{eq: smoothness_2}} \text{ holds},\forall x_1,x_2\in \Omega\}$. 
\end{assumption}
In the remainder of this paper, we also define $L_{\max} = \max_{i\geq 1}L_i$ and  $M_{\max} = \max_{i\geq 1}M_i$. 

\begin{assumption}
\label{ass: strict feasible}
There exists a known strictly feasible point $x_0$, i.e., $f_i(x_0)<0$ for all $i\in\mathbb Z_1^m$. 
\end{assumption}

\begin{assumption}
\label{ass: bounded_sublevel}
There exists $\beta\in \mathbb{R}$ such that the sublevel set $\mathcal{P}_\beta=\{x\in \Omega:f_0(x)\leq \beta\}$ is bounded and includes the initial feasible point $x_0$.
\end{assumption}

Assumption \ref{ass: strict feasible} is common in safe zeroth-order methods \cite{sui2015safe,usmanova2019safe,usmanova2022log,guo2023safe}. Without the initial feasible point, it would be impossible to ensure the feasibility of all the samples. Assumption \ref{ass: bounded_sublevel} is not strong since it holds as long as the feasible region $\Omega$ is bounded. 

Our aim is to derive an optimization algorithm where a subsequence of the iterates converges to the primal of a KKT pair.
\begin{definition}
\label{def: KKT}
If a pair $(x,\lambda)$ with $x\in \Omega$ and $\lambda\in \mathbb{R}^m_{\geq 0}$ satisfies 
\begin{subequations}
\label{eq:kkt}
\begin{align}
\|\nabla f_0(x) + \sum^m_{i=1}\lambda^{(i)}\nabla f_i(x)\| &=0,\\
|\lambda^{(i)}f_i(x)| &= 0,\quad i\in \mathbb{Z}^{m}_{1},
\end{align}
\end{subequations}
we say that $(x,\lambda)$ is a KKT pair of the problem \eqref{eq: optimization problem} and $x\in\Omega$ is the primal of the KKT pair. 
\end{definition}

In the following section, we design a safe zeroth-order algorithm whose iterates, under mild assumptions, have an accumulation point that is the primal of a KKT pair of \eqref{eq: optimization problem}.

\section{Algorithm: SZO-LP}
\label{sec: algorithm}
In this section, we first describe how to estimate gradients of the functions in \eqref{eq: optimization problem} and construct local feasible sets. These are the essential tools used by our zeroth-order optimization method.
\subsection{Gradient estimation and local feasible set construction}
We estimate the gradient through finite difference, i.e., 
\begin{equation}
\label{eq: gradient approximation}
{\nabla}^\nu f_i\left(x\right):=\sum_{j=1}^{d} \frac{f_i\left(x+\nu e_{j}\right)-f_i\left(x\right)}{\nu} e_{j}.
\end{equation}
The following lemma gives a method to control the estimation error $\Delta^\nu_i(x):={\nabla}^\nu f_i\left(x\right)-\nabla f_i(x).$
\begin{lemma}[\cite{berahas2021theoretical}, Theorem 3.2]
\label{lmm: gradient estimation error}
Under Assumption \ref{ass: smoothness}, we have
\begin{equation}
\left\|\Delta^\nu_i(x)\right\|_{2} \leq \frac{\sqrt{d}  M_i}{2} \nu.
\label{eq: gradient approximation error}
\end{equation}
By letting $\nu = \nu(\epsilon):= \frac{2\epsilon}{\sqrt{d}M_{\max}}$ we have
$\|\Delta^{\nu(\epsilon)}_i(x)\|\leq \epsilon$.
\end{lemma}

By using the estimated gradient, we can build a local feasible set around $x_0$. 
We let \begin{equation}
\label{eq: l_0}
l_0^* = \min_{i\in\{1, \dots, m\}}\;-f_i(x_0)/L_{\mathrm{max}}, 
\end{equation}
and $\nu_0^*(\epsilon) := \min\{l_0^*/\sqrt{d},\nu(\epsilon)\}$. Since $\nu_0^*(\epsilon)\leq \nu(\epsilon)$, we have $\|\Delta^{\nu_0^*(\epsilon)}_i(x_0)\|\leq \epsilon$ for any $\epsilon>0$ and $i\geq 1$.

Then the set
\begin{equation}
\label{eq:safe_set}
\begin{aligned}
\mathcal{S}^{(0)}(x_0) := & \cap^{m}_{i=1}\mathcal{S}^{(0)}_i(x_0), \text{ where}\\
\mathcal{S}^{(0)}_i(x_0) :=&\big\{x: f_i(x_0)+{\nabla}^{\nu_0^*(\epsilon_0)} f_i\left(x_0\right)^\top(x-x_0)+\\
&\qquad\qquad\qquad\qquad 2M_i\|x-x_0\|^2\leq 0 \big\}
\end{aligned}
\end{equation}
is feasible as shown in the following theorem.
\begin{theorem}[\cite{guo2023safe}, Theorem 1]
\label{thm: local feasible set}
All the samples used to construct $\mathcal{S}^{0}(x_0)$ are feasible. Moreover, the set $\mathcal{S}^{(0)}(x_0)$ is convex and any $x\in \mathcal{S}^{(0)}(x_0)$ is strictly feasible. 
\end{theorem}
In the lack of explicit constraint functions, a local feasible set is a common tool of several zeroth-order methods \cite{sui2015safe,usmanova2022log,guo2023safe} to ensure the feasibility of the iterates, though the specific formulations are different. In the following, we propose our method where the local feasible sets are used to select the step length for the derived descent direction. 

\subsection{Algorithm: Safe Zeroth-Order Optimization Using Linear Programs (SZO-LP)}
The main idea of the SZO-LP method, shown in Algorithm~\ref{alg: SZO-LP}, is to iteratively select a descent direction by executing in Line 7 $\mathtt{LP}(x_k,\epsilon_k)$ defined in \eqref{eq: LP_definition}. Thanks to the tightening contant $\epsilon_k$ in the linear program involved in $\mathtt{LP}(x_k,\epsilon_k)$, the descent direction we obtain points into the iterior of the feasible set. Along this direction, we select the step length (Line 9-14) based on local feasible sets and the pre-defined length
$$\gamma(\epsilon_k) := \frac{\epsilon_k}{4(M_{\max}+L_{\max})}.$$
\begin{algorithm}[htbp!]
\caption{Safe Zeroth-Order optimization using Linear Programs (SZO-LP)}
\label{alg: SZO-LP}
\textbf{Input:} $\epsilon_0$, $\epsilon_{\min}$, $K_{\mathrm{switch}}$, initial feasible point $x_0\in\Omega$\\
\textbf{Output:} $\tilde{x}$
\begin{algorithmic}[1]
\State $k\gets 0, \textsc{TER}= 0$ 
\While{$\epsilon_k>\epsilon_{\min}$}
\State $s_\mathrm{tmp} \gets \mathtt{LP}(x_k,2\epsilon_k)$
\If{$\nabla^{\nu^*_k(2\epsilon_k)}f_0(x_{k})^\top s_\mathrm{tmp}\leq -4\epsilon_k$}
\State $\epsilon_{k+1}\gets 2\epsilon_k$, $x_{k+1}\gets x_{k}$
\Else
\State $s_k^* = \mathtt{LP}(x_k,\epsilon_k)$
\vspace{0.1cm}
\If{$\nabla^{\nu^*_k(\epsilon_k)}f_0(x_{k})^\top s^*_k\leq -2\epsilon_k$ } 
\vspace{0.1cm}
\If{$k<K_{\mathrm{switch}}$}
\begin{small}
\begin{align}
\label{eq: step_length_1}
&\beta_k =\argmax_{\beta\geq0} \beta \text{ s.t. }x_k+\beta s^*_k\in \mathcal{S}^{(k)}(x_k), \\
\label{eq: step_length}
&\alpha_k = \argmin_{\alpha\in \{\beta_k,\gamma(\epsilon_k)\}} f_0(x_k+\alpha s_k^*)
\end{align}
\end{small}
\State $x_{k+1}\gets x_k+\alpha_k s^*_k$, $\epsilon_{k+1}\gets \epsilon_{k}$
\Else
\State $x_{k+1}\gets x_k+\gamma(\epsilon_k) s^*_k$, $\epsilon_{k+1}\gets \epsilon_{k}$
\EndIf
\Else 
\State $\epsilon_{k+1}\gets \epsilon_k/2$, $x_{k+1}\gets x_{k}$
\EndIf

\EndIf

\State $k\gets k+1$
\EndWhile
\end{algorithmic}
\end{algorithm}

The essential steps are as follows:
\subsubsection{Providing the input data}
The input includes an initial strictly feasible point $x_0$ (see Assumption \ref{ass: strict feasible}) and a tightening constant $\epsilon_0$. Each iteration of the algorithm generates a new tightening constant $\epsilon_k$, which can be equal to $\epsilon_{k-1}$, $2\epsilon_{k-1}$ or $\epsilon_{k-1}/2$. Since $\epsilon_k$ converges to 0 (see Theorem \ref{thm: tools_for_convergence}), the user can control the termination by providing a lower bound $\epsilon_{\min}$ for $\epsilon_k$. The parameter $K_{\text{switch}}$ marks the boundary of two methods for selecting step length, see the last bullet point.
\subsubsection{Building local feasible sets} 
For a strictly feasible $x_k$, we use \eqref{eq: l_0} to define $l_k^*$ and 
$$\nu_k^*(\epsilon_k) := \min\{l_k^*/\sqrt{d},\nu(\epsilon_k)\}.$$ We then use $\nu_k^*(\epsilon_k)$ and \eqref{eq:safe_set} to define
$\mathcal{S}^{(k)}(x_k)$, a local feasible set around $x_k$. From Theorem \ref{thm: local feasible set} we know that if $x_{k+1}\in \mathcal{S}^{(k)}(x_k)$ then $x_{k+1}$ is also strictly feasible.
\subsubsection{Solving subproblems for the descent diretion} 
In each iteration, we execute in Line 7 \texttt{LP}$(x_k,\epsilon_k)$ to derive a search direction, which returns
\begin{equation}
\label{eq: LP_definition}
\begin{aligned}
\argmin_{\|s\|_1\leq 1}\;\;&(\nabla^{\nu^*_k(\epsilon_k)}f_0(x_{k}))^\top s\\
\mathrm{s.t.}\;\;&(\nabla^{\nu^*_k(\epsilon_k)}f_i(x_{k}))^\top s+2\epsilon_k \leq 0,\\
 &\, \forall i \in \mathcal{A}(x_k,\epsilon_k),
\end{aligned}
\end{equation}
or NaN if \eqref{eq: LP_definition} is not feasible. Here, $\mathcal{A}(x,\epsilon):=\{i:f_i(x)\geq -2\epsilon\}$ is the near-active constraint index set. The solution to \eqref{eq: LP_definition} is a direction that not only gives a fast descent but also points into the interior of the feasible region $\Omega$ (away from the boundary). In \eqref{eq: LP_definition}, due to the tightening constant $\epsilon_k$, along the direction $s_k^*$ in Line 7, the constraint function values decrease. Therefore, moving along the direction $s_k^*$ we indeed stay away from the boundary of $\Omega$. This direction helps to avoid small values of $-f_i(x_k)$, which lead to conservative local feasible sets $\mathcal{S}^{(k)}(x_k)$. Moreover, the inclusion of only near-active constraints makes \eqref{eq: LP_definition} small-size and easy to solve. We will later see in Theorem \ref{thm: tools_for_convergence} that $\epsilon_k$ converges to 0. Therefore, it is still possible that a subsequence of the iterates converges to a point on the feasible set boundary.  

We also let $s_\mathrm{tmp} = \mathtt{LP}(x_k,2\epsilon_k)$ and check in Line~4 whether $\nabla^{\nu^*_k(2\epsilon_k)}f_0(x_{k})^\top s_\mathrm{tmp}\leq -4\epsilon_k$, which allows us to have Proposition \ref{prop: lower_bound_for_epsilon}, the proof of which is in Appendix \ref{sec: proof of lower_bound_for_epsilon}. This proposition will be later used to show in Theorem \ref{thm: kkt} the properties of the $\{x_k\}_{k\geq 1}$ as $k$ goes to infinity.
\vspace{-0.3cm}
\begin{proposition}
\label{prop: lower_bound_for_epsilon}
Any $\epsilon_k$ entering Line 7 satisfies 
\begin{align}\label{eq: largest epsilon}
\epsilon_k \geq \frac{1}{8}{\mathrm{sup}}\{\epsilon: & s = \mathtt{LP}(x_k,\epsilon) \textrm{ verifies} \\\notag
&\qquad\qquad\qquad\nabla^{\nu^*_k(\epsilon)}f_0(x_{k})^\top s\leq -2\epsilon\}.  
\end{align}
\end{proposition}

\subsubsection{Deciding the step length} When a direction $s^*_k$ derived in Line 7 gives sufficient descent (i.e., $\nabla^{\nu^*_k(\epsilon_k)}f_0(x_{k})^\top s^*_k\leq  -2\epsilon_k$), we move along the tentative direction $s^*_k$. To decide the step length, we consider the local feasible set and the pre-defined step length $\gamma(\epsilon_k)$ that is guaranteed to achieve a non-trivial descent (see Lemma \ref{lmm: sufficient_decrease}). In \eqref{eq: step_length_1}, we calculate by bisection the largest step length within the local feasible set to derive $\alpha_k$ in \eqref{eq: step_length}. The use of local feasible sets in Line 10 allows us to obtain a larger step length than $\gamma(\epsilon)$, when $x_k$ is not close to the boundary of the feasible set. This is because, from the formulation \eqref{eq:safe_set}, smaller values of $f_i(x_k)$ lead to larger sizes of $\mathcal{S}^{(k)}_i(x_k)$ while $\gamma(\epsilon_k)$ is independent of how far the iterates are from the feasible set boundary. When $k>K_{\mathrm{switch}}$, we let the step length be $\gamma(\epsilon_k)$ as in Line 12, which is useful for the proof of the iterates' properties as $k$ goes to infinity (see Theorem \ref{thm: kkt}). The selection of $K_{\mathrm{switch}}$ is not critical since we use the step length in Line 10 for $k<K_{\mathrm{switch}}$ instead of that defined in Line 12 only to accelerate the descent in the early iterations of the algorithm. 

On the other hand, if the direction $s^*_k$ cannot give sufficient descent, we let $\epsilon_{k+1} = \epsilon_k/2$ in Line 15 to relax the tightened constraints in \eqref{eq: LP_definition}. This relaxation makes it easier for $s^*_{k+1}$ to give sufficient descent, i.e., to satisfy $\nabla^{\nu^*_k(\epsilon_{k+1})}f_0(x_{k+1})^\top s^*_{k+1}\leq  -2\epsilon_{k+1}$. Only when $s^*_{k+1}$ gives sufficient descent will we move along $s^*_{k+1}$ to a new point.
\vspace{0.2cm}

We refer the readers to Remark \ref{rmk: comparison} for how SZO-LP is compared with some state-of-the-art methods.

\section{Convergence properties of the approach}
\label{sec: theories}
In this section, we aim to show that, under mild conditions and by letting $\epsilon_{\min}=0$, the sequence $\{x_k\}_{k\geq1}$ produced in Algorithm \ref{alg: SZO-LP} has an accumulation point $x_c$ that is also the primal of a KKT pair of \eqref{eq: optimization problem}. To start with, we show in Lemma \ref{lmm: sufficient_decrease} that, whenever $x_{k+1}\neq x_{k}$, the new iterate $x_{k+1}$ is strictly feasible and the objective function value gets a non-trivial decrease.
\begin{lemma}
\label{lmm: sufficient_decrease}
Suppose $s_k^*$ derived in Line 7 of Algorithm \ref{alg: SZO-LP}, satisfies $$\nabla^{\nu^*_k(\epsilon_k)}f_0(x_{k})^\top s_k^*\leq -2 \epsilon_k.$$ We have that $x_k+\gamma(\epsilon_k)s_k^*$ is strictly feasible. Furthermore $x_k+\gamma(\epsilon_k)s_k^*$ satisfies
\begin{equation}
f_0(x_k+\gamma(\epsilon_k)s_k^*)-f_0(x_k)< -\epsilon^2_k/(8(M_{\max}+L_{\max})). 
\end{equation}
\end{lemma}
The proof of Lemma \ref{lmm: sufficient_decrease} is in Appendix \ref{sec: sufficient_decrease}. The main idea is to utilize the smoothness constants in Assumption \ref{alg: SZO-LP} to upper-bound $f_i(x_k+\gamma(\epsilon_k))$ for $i\in\mathbb{Z}^m_0 $. Based on this lemma, we have the following theorem on the sequences $\{x_k\}_{k\geq 1}$ and $\{\epsilon_k\}_{k\geq 1}$ as $k$ goes to infinity. 

\begin{theorem}
\label{thm: tools_for_convergence}
The following arguments hold:
\begin{itemize}
    \item[1.] The sequence $\{f_0(x_k)\}_{k\geq 1}$ is non-increasing;
    \item[2.] There exists at least one accumulation  point of the sequence $\{x_k\}_{k\geq 1}$. For any accumulation point $x_c$,
$$\lim_{k\rightarrow \infty} f_0(x_k) = f_0(x_c)>-\infty.$$
    \item[3.] The sequence $\{\epsilon_k\}_{k\geq 1}$ converges to 0.
\end{itemize}
\end{theorem}
\textbf{Proof.} The first point is a direct consequence of Lemma \ref{lmm: sufficient_decrease}, which implies that whenever the iterate moves to a new point the objective function value decreases. 

\textit{Proof of Point 2.} Since $\{f_0(x_k)\}_{k\geq 1}$ is non-increasing, $f_0(x_k)\leq f_0(x_0)$, for any $k\geq 1$, and thus $x_k\in \mathcal{P}_\beta$. Due to the boundedness of $\mathcal{P}_\beta$, by the Bolzano–Weierstrass theorem, we know that there exists at least one accumulation point of $\{x_k\}_{k\geq 1}$. For any accumulation point $x_c$, there exists a subsequence $\{x_{k_p}\}_{p\geq 1}$ converging to $x_c$. Due to the continuity of $f_0(x)$, 
$$\lim_{p\rightarrow \infty} f_0(x_{k_p}) = f_0(x_c)>-\infty.$$
By utilizing again the monotonicity of $\{f_0(x_k)\}_{k\geq 1}$, we have 
$\lim_{k\rightarrow \infty} f_0(x_k) = f_0(x_c)>-\infty. $ 

\textit{Proof of Point 3.} We show this result through contradiction by assuming that $\{\epsilon_k\}_{k\geq 1}$ does not diminish as $k$ goes to infinity. Based on this assumption, we can show that $\{\epsilon_k\}_{k\geq 1}$ does not converge to any non-zero values. If $\{\epsilon_k\}_{k\geq 1}$ converges to a non-zero value, by noticing that $\epsilon_k\in\{\epsilon_0*2^i:i\in\mathbb{Z}\}$,  we have that there exists $K>0$ such that any $k>K$ verifies $\epsilon_k = \epsilon_{k-1}$. Then for any $k>K$ the new iterate $x_{k+1}$ is derived in Line 10 or 12 in Algorithm \ref{alg: SZO-LP}. According to Lemma \ref{lmm: sufficient_decrease}, $f_0(x_{k+1})\leq f_0(x_k)-\epsilon^2_k/(8(M_{\max}+L_{\max}))$ and thus $f_0(x_k)$ goes to $-\infty$ as $k$ goes to $+\infty$, which contradicts Point 2. Therefore, $\epsilon_{k}$ does not converge. 

Since from Algorithm \ref{alg: SZO-LP} $\{\epsilon_k\}_{k\geq 1}$ is bounded, we can conclude that $\{\epsilon_k\}_{k\geq 1}$ has multiple accumulation points. Then there are $\epsilon>0$ and infinitely many $k$ such that $\epsilon_k = \epsilon$ and $\epsilon_{k-1} = \epsilon/2$. For any $k$ of this kind, there exists $k'\geq k$ verifying $f_0(x_{k'+1})\leq f_0(x_{k'})-\epsilon^2/(8(M_{\max}+L_{\max}))$. Consequently there are infinitely many $k'$ verifying $f_0(x_{k'+1})\leq f_0(x_{k'})-\epsilon^2/(8(M_{\max}+L_{\max}))$, which again contradicts Point 2. \hfill $\blacksquare$
\vspace{0.3cm}

Theorem \ref{thm: tools_for_convergence} offers us the essential tools to show in Theorem~\ref{thm: kkt} the properties of an accumulation point of $\{x_k\}_{k\geq 1}$ under Assumption \ref{ass: LICQ}.

\begin{assumption}
\label{ass: LICQ}
At least one accumulation point $x_c$ of $\{x_k\}_{k\geq 1}$ satisfies \textit{Linear Independent Constraint Qualification (LICQ)}, which is to say the gradients $\nabla f_i(x_c)$ with $i\in \mathcal{A}(x_c,0)$ are linearly independent. 
\end{assumption}
Assumption \ref{ass: LICQ} is widely used in optimization \cite{wachsmuth2013licq}. For example, it is used to prove the properties of the limit point of the Interior Point Method \cite{jorge2006numerical}.

\begin{theorem}
\label{thm: kkt}
Regarding the accumulation point $x_c$ in Assumption \ref{ass: LICQ}, there exists $\lambda_c\in\mathbb{R}^m_{\geq 0}$ such that $(x_c,\lambda_c)$ is a KKT pair of \eqref{eq: optimization problem}.
\end{theorem}
The proof, in Appendix \ref{sec: kkt}, is based on contradiction. If $x_c$ is not the primal of a KKT pair, we can find $r>0$, $\epsilon>0$ and $s_\epsilon\in \mathbb{R}^d$  such that for any $x_k\in\mathcal{B}_r(x_c)$ the solution $s = \mathtt{LP}(x_k,\epsilon) $ verifies $\nabla^{\nu^*_k(\epsilon)}f_0(x_{k})^\top s\leq -2\epsilon$. There are infinitely many $k$ such that $x_k\in\mathcal{B}_r(x_c)$ and $s^*_{k}$ is derived through Line 7 in Algorithm \ref{alg: SZO-LP}. For any of these $k$s, according to \eqref{eq: largest epsilon}, $\epsilon_k\geq \epsilon/8$, which contradicts Point 3 of Theorem \ref{thm: tools_for_convergence}.
\vspace{0.3cm}

\begin{remark}
\label{rmk: comparison}
Like SZO-QQ \cite{guo2023safe} and LB-SGD \cite{usmanova2022log}, the samples in SZO-LP are all feasible and the iterates, under mild assumptions, have an accumulation point that is also the primal of a KKT pair. In contrast, the tightening constant $\epsilon_k$ of SZO-LP keeps the iterates away from the boundary of the feasible set and leads to less conservative local feasible sets than those used in SZO-QQ and LB-SGD. Moreover, due to the use of the near-active set $\mathcal{A}(x_k,\epsilon_k)$ the subproblems \eqref{eq: LP_definition} are smaller-size and easier to solve than the QCQPs in SZO-QQ and nonconvex subproblems in Safe Bayesian Optimization methods \cite{sui2015safe,berkenkamp2021bayesian}. However, to rigorously show these advantages, we need to upper bound the number of iterations needed by SZO-LP given certain accuracy requirements, which is left as future work.
\end{remark}

\section{Experiment on an OPF problem}
\label{sec: numericals}
To illustrate the performance of SZO-LP, we consider applying it to an OPF problem on the IEEE 30-bus system.
\subsection{Formulation of the OPF problem}
To formulate an OPF problem, we introduce the following notations and assumptions:
\begin{itemize}
    \item Let $B=\{b_1,b_2,\ldots,b_n\}$ be the bus set and let $T=\{(b_i,b_j):$ there is a transmission line between $b_i$ and $b_j\}$ be a set of undirected edges representing the transmission lines;
    \item We denote $P_{G_i}$, $P_{L_i}$, $Q_{L_i}$, $U_i$ and $\theta_i$ as the active power generation, active power consumption, reactive power consumption, voltage and voltage angle at $b_i$;
    \item From the $b_i$ to $b_j$, the active power and the reactive power transferred are written respectively as $P_{ij}(U_i,U_j,$ $\theta_i,\theta_j)$ and $Q_{ij}(U_i,U_j,\theta_i,\theta_j)$, while the current is denoted as $I_{ij}(U_i,U_j,\theta_i,\theta_j)$. We refer the readers to \cite{das2017load} for the explicit expressions of these functions;
    \item We also assume that there are $n_G$ generators at the buses $b_i$, $i\in \mathbb{Z}^{n_G}_1$ and $b_1$ is a slack bus providing active power to maintain the power balance within the network and has a voltage angle of 0.
\end{itemize}
Then the OPF problem is formulated \cite{das2017load} as 
\begin{subequations}
\label{eq:OPF}
\begin{align}
\min_{P_{G_i},U_i,\theta_i} & \sum_{i=1}^{n_G}C_i(P_{G_i})
\end{align}
subject to
\begin{align}
 P_{G_i} &= P_{L_i}+\sum_{(i,j)\in T}P_{ij}(U_i,U_j,\theta_i,\theta_j),\;\; \forall i\label{eq: power_flow_active}\\
 -Q_{L_i} &= \sum_{(i,j)\in T}Q_{ij}(U_i,U_j,\theta_i,\theta_j),\;\; i>n_G \label{eq: power_flow_reactive}\\
 P_{G_i}& = 0,\;\; \mathrm{for }\text{ } i>n_G, \;\;\theta_1 = 0, \label{eq: power_flow_generation}\\
 P_{G,\min}&\leq P_{G_i} \leq P_{G,\max}, \mathrm{for} \;\;i\leq n_G, \label{eq: P_range}\\
 I_{ij,\min}&\leq  I_{ij}(U_i,U_j,\theta_i,\theta_j)\leq I_{ij,\max},\;\; \forall (i,j)\in T ,\label{eq: I_range}\\
 U_{\min}&\leq U_{i} \leq U_{\max}, \forall i. \label{eq: U_range}
\end{align}
\end{subequations}
where $C_i(\cdot)$ is a quadratic function accounting for the generation cost and the equations \eqref{eq: P_range}-\eqref{eq: U_range} give the safe intervals for the corresponding variables. 

The main challenges of OPF applications lie in modelling the system and deriving the accurate expressions of \eqref{eq:OPF}. The difficulties include the nonlinearity of device dynamics, slowly changing physical parameters and disturbances \cite{chuMitigating22}. Inaccurate models can result in suboptimal OPF solutions (leading to more generation cost) or violate the true hard constraints (causing damages to devices) \cite{lee2021robust}. Therefore, we consider the black-box setting and use SZO-LP. 

To this aim, we reformulate \eqref{eq:OPF} as optimization with only inequality constraint to fit \eqref{eq: optimization problem} used by SZO-LP. Let $\{P_{G_i}\}^{n_G}_{i=2}$ and $\{U_i\}^{n_G}_{i=1}$ be the main decision variables. Then by assigning values to $\{P_{G_i}\}^{n_G}_{i=2}$ and $\{U_i\}^{n_G}_{i=1}$, one can solve the power flow equations \eqref{eq: power_flow_active}-\eqref{eq: power_flow_generation} to derive the values for all the other decision variables in \eqref{eq:OPF}. Therefore, \eqref{eq: power_flow_active}-\eqref{eq: power_flow_generation} give us the functions
    \begin{equation}
    \label{eq: implicit functions}
    \begin{aligned}
    U_{i} &= U_{i}(\{P_{G_j}\}^{n_G}_{j=2},\{U_j\}^{n_G}_{j=1}),\;\;i=n_G+1,\ldots, n, \\
    \theta_{i} &= \theta_{i}(\{P_{G_j}\}^{n_G}_{j=2},\{U_j\}^{n_G}_{j=1}),\;\;i=1,\ldots, n.
        \end{aligned}
    \end{equation}
By substituting \eqref{eq: implicit functions} to \eqref{eq:OPF}, we obtain a reformulation where $\{\{P_{G_j}\}^{n_G}_{j=2},\{U_j\}^{n_G}_{j=1}\}$ are the only decision variables and there are not equality constraints. 

\subsection{Experiment results}
We run SZO-LP to solve a specific OPF problem on the IEEE 30-bus system where $n_G = 6$ . In total, there are 11 decision variables and 158 constraints. We do not assume knowledge of the system model for the optimization task. However, given a set of values for all 11 decision variables, we can use a black-box simulation model in \texttt{Matpower}~\cite{zimmerman2010matpower} to sample the voltages of all the 30 buses and the current through all the transmission lines in the network. Additionally, we assume the availability of initial values for all the decision variables to start the SZO-LP algorithm from a feasible point.

We employ SZO-LP to reduce the quadratic cost induced by the initial decision values. The numerical experiments are executed on a PC with an Intel Core i9 processor. The solver we adopt for subproblems \eqref{eq: LP_definition} is \textsc{linprog} in Matlab. We let $M_i = M_{\max}=0.13$ and $L_i = L_{\max}=0.5$. The tuning of these two parameters is described in \cite{guo2023safe}. Moreover, we set $\epsilon_0=0.05$, $\epsilon_{\min} = 10^{-6}$ and $K_{\mathrm{switch}} = 200$. 


\begin{figure}[htbp!]
\centering
\includegraphics[width=0.8\linewidth]{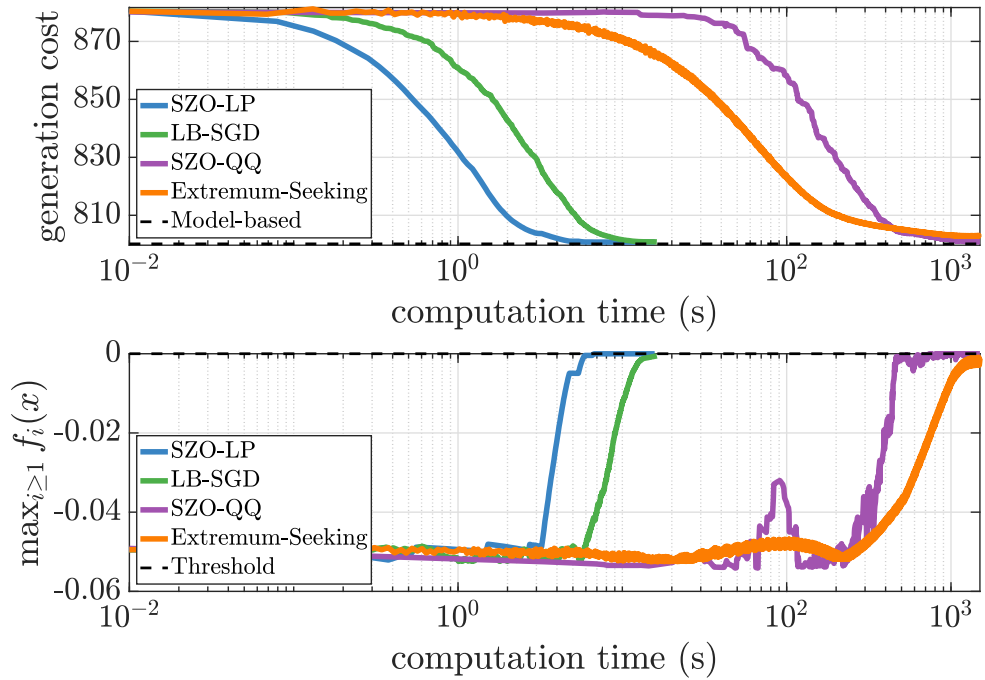}
\caption{Decrease of cost and growth of the largest constraint function values with respect to computation time}
\label{fig:opf_comparison_time}
\end{figure}



\begin{figure}[htbp!]
\centering
\includegraphics[width=0.8\linewidth]{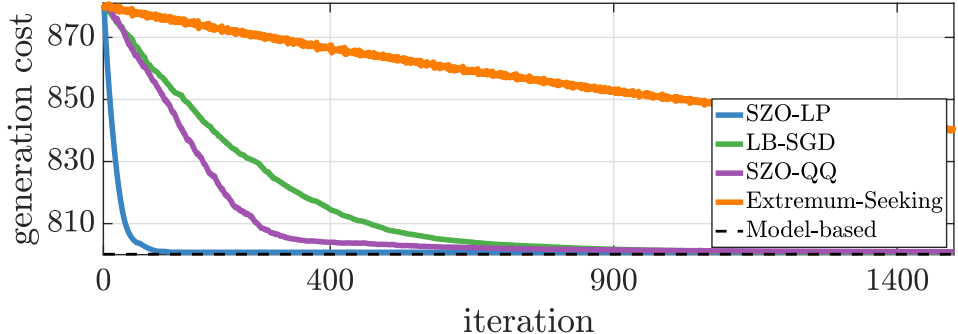}
\caption{Decrease of cost with respect to the number of iterations}
\label{fig:opf_iteration}
\end{figure}

In Figures \ref{fig:opf_comparison_time} and  \ref{fig:opf_iteration} , we present the results of our numerical experiments, where we compare the performance of SZO-LP with SZO-QQ \cite{guo2023safe}, LB-SGD \cite{usmanova2022log} and Extremum Seeking \cite{arnold2015model}. The QCQP subproblems in SZO-QQ are solved using MOSEK. The reference solution of the OPF problem is returned by the optimization based on the true model and utilizing Gurobi~\cite{gurobi} as the solver. The computation time in Figure \ref{fig:opf_comparison_time} includes that consumed by power grid simulation (through Matpower) when we query the objective and constraint functions. We observe that all four methods keep the iterates feasible and eventually achieve a generation cost very close to that (800.14) derived based on the true model. However, SZO-LP achieves a faster decrease in the generation cost than the other methods.

One main reason for the superior performance of SZO-LP over SZO-QQ with respect to computation time shown in Figure \ref{fig:opf_comparison_time},  is that the linear programming subproblems can be solved faster. We notice that to finish the first 60 subproblems, SZO-LP takes 5.63 seconds while SZO-QQ takes 72.06 seconds. Firstly, the subproblem in SZO-LP only takes into account the near-active constraints while the subproblem in SZO-QQ involves all constraints. Among the iterations of SZO-LP, the largest number of constraints is 2. Secondly, although the big gap in efficiency shown in Figure \ref{fig:opf_comparison_time} may be due to the specific solvers we select, linear programs, in general, are open to a wider selection of solvers and thus allow for more efficient implementations. 

Unlike SZO-LP and SZO-QQ, LB-SGD and Extremum Seeking do not require solving any subproblems, thus allowing for more iterations within a certain time length. This is why LB-SGD can also achieve a low generation cost in a short time. However, considering the four methods take the same number of samples every iteration, LB-SGD and Extremum Seeking are less sample-efficient than SZO-LP and SZO-QQ since they require more iterations as shown in Figure \ref{fig:opf_iteration}. Moreover, since LB-SGD and Extremum Seeking are based on log barriers, these two methods require tuning of the barrier function coefficients. Improper tuning might lead to suboptimality in LB-SGD or even infeasibility in Extremum Seeking.

SZO-LP has another advantage over SZO-QQ, which is the feature of SZO-LP keeping the iterates away from the feasible set boundary before getting close to the primal of a KKT pair. Iterates getting too close to the feasible set boundary might impede the decrease of the cost. To see this point, we notice from Figure 
\ref{fig:opf_comparison_time} that in SZO-QQ the decrease of the generation cost slows down when the largest constraint function value is larger than -0.005. The reason is that, when the largest constraint function value is close to 0, the local feasible set constructed in SZO-QQ gets conservative, and thus the step length becomes small. When the largest constraint function value gets larger than -0.005 for the first time, the generation cost in SZO-QQ is 805.27 while the corresponding cost in SZO-LP is 801.77, which is much closer to 800.14 (derived by optimization based on the true model). Therefore, we see that in SZO-QQ the decrease of the objective function value can slow down at a much earlier stage. 

In conclusion, from the experiment results, we see that SZO-LP is the most computation-efficient and sample-efficient method, among the four approaches.

\section{Conclusion}
In this paper, we proposed a safe zeroth-order method SZO-LP, which iteratively solves linear programs to obtain descent directions and determines the step lengths. We showed that, under mild conditions, the iterates of SZO-LP have an accumulation point that is also the primal of a KKT pair. Through an experiment where we use SZO-LP to solve an OPF problem on the IEEE 30-bus system and compare with three other methods, we see that SZO-LP is both computation-efficient and sample-efficient. Our future directions include the derivation of the computation complexity of SZO-LP to check whether it is efficient in general and the extension of SZO-LP to account for measurement noises.

\bibliographystyle{ieeetr}
\bibliography{reference}

\section*{Appendix}
\appendix
{
\section{The proof of Proposition \ref{prop: lower_bound_for_epsilon}}
\label{sec: proof of lower_bound_for_epsilon}
We first show that if for some $k>0$ and $\epsilon_\alpha>0$
\begin{equation}
\label{eq: inequality_1}
s_1 = \mathtt{LP}(x_k,\epsilon_\alpha) \textrm{ verifies } \nabla^{\nu^*_k(\epsilon_\alpha)}f_0(x_{k})^\top s_2\leq -2\epsilon_\alpha    
\end{equation} 
then for any $\epsilon_\beta\leq \epsilon_\alpha/4$
\begin{equation}
\label{eq: inequality_2}
s_2 = \mathtt{LP}(x_k,\epsilon_\beta) \textrm{ verifies } \nabla^{\nu^*_k(\epsilon_\beta)}f_0(x_{k})^\top s_2\leq -2\epsilon_\beta.
\end{equation} 
With \eqref{eq: inequality_1}, we notice that $s_1$ with $\|s_1\|_1\leq 1$ is a feasible solution to the linear program involved in $\mathtt{LP}(x_k,\epsilon_\beta)$. This is because for any $i\in \mathcal{A}(x_k,\epsilon_\beta)\subset \mathcal{A}(x_k,\epsilon_\alpha)$, we have
\begin{equation*}
\begin{aligned}
&\langle \nabla^{\nu^*_k(\epsilon_\beta)} f_i(x_k), s_1 \rangle\\
\leq & \langle \nabla f_i(x_k), s_1 \rangle + |\langle\nabla f_i(x_k)-\nabla^{\nu^*_k(\epsilon_\beta)} f_i(x_k), s_1 \rangle|\\
\leq & \langle \nabla^{\nu^*_k(\epsilon_\alpha)} f_i(x_k), s_1 \rangle + |\Delta_i^{\nu_k^*(\epsilon_\alpha)}(x_k)| + |\Delta_i^{\nu_k^*(\epsilon_\beta)}(x_k)|\\
\leq & -2\epsilon_\alpha+\epsilon_\alpha+\epsilon_\alpha/4 <-\epsilon_\alpha/2\leq -2\epsilon_\beta.
\end{aligned}
\end{equation*}
Similarly, we can show that 
$$\langle \nabla^{\nu^*_k(\epsilon_\beta)} f_0(x_k), s_1 \rangle \leq -2\epsilon_\beta.$$
Considering that $s_2$ is the optimum of the linear program involved in $\mathtt{LP}(x_k,\epsilon_\beta)$, \eqref{eq: inequality_2} holds. 

Then if $\epsilon_k$ enters Line 7 of Algorithm \ref{alg: SZO-LP}, the condition in Line 4 
``$s_\mathrm{tmp} = \mathtt{LP}(x_k,2\epsilon_k)$ verifying $\nabla^{\nu^*_k(2\epsilon_k)}f_0(x_{k})^\top s_\mathrm{tmp}\leq -4\epsilon_k$'' does not hold. By letting 
$$\epsilon_\alpha = \mathrm{sup}\{\epsilon:  s = \mathtt{LP}(x_k,\epsilon) \textrm{ verifies } \nabla^{\nu^*_k(\epsilon)}f_0(x_{k})^\top s\leq -2\epsilon\},$$
for any $\epsilon_\beta\leq \epsilon_\alpha/4$, \eqref{eq: inequality_2} holds. Therefore, $2\epsilon_k\geq \epsilon_\alpha/4$. \hfill $\blacksquare$
\section{Proof of {Lemma~\ref{lmm: sufficient_decrease}}}
\label{sec: sufficient_decrease}
To begin with, we show that $x_k+\gamma(\epsilon_k)s_k^*$ is indeed strictly feasible. By using the mean value theorem and noticing that $\|s_k^*\|\leq \|s_k^*\|_1\leq 1$, we have that for any $\gamma>0$
\begin{align}\notag
& f_i(x_k+\gamma s_k^*)\\\notag
< & f_i(x_k)+\gamma \nabla f_i(x_{k})^\top s_k^* + 2\gamma^2M_{\max}\|s_k^*\|^2 \\\notag
< &  \gamma \nabla^{\nu^*_k(\epsilon_k)} f_i(x_{k})^\top s_k^* + \gamma \|\Delta^{\nu^*_k(\epsilon_k
)}_i(x)\| \cdot \|s_k^*\| + 2M_{\max}\gamma^2 \\\notag
< & -2\epsilon_k \gamma + \epsilon_k \gamma  + 2 M_{\max}\gamma^2\\\label{eq: upperbounding f_i}
< &  2 (M_{\max}+L_{\max}) \gamma ^2 -\epsilon_k \gamma, \text{ }\forall i\in \mathcal{A}(x_k,\epsilon_k),\\
\notag&f_i(x_k+\gamma s_k^*) \\
\label{eq: upperbounding f_i_non_active}<&f_i(x_k)+L_{\max}\gamma, \text{ }\forall i\in \mathbb{Z}^m_1\setminus\mathcal{A}(x_k,\epsilon_k).
\end{align}
Therefore, we have $$f_i(x_k+\gamma(\epsilon_k)s_k^*)<-\epsilon^2_k/(8(M_{\max}+L_{\max}))<0$$ for any $i\in \mathcal{A}(x_k,\epsilon_k)$ and 
$$f_i(x_k+\gamma(\epsilon_k)s_k^*)<-\epsilon_k/2$$ 
for any $i\in \mathbb{Z}^m_1\setminus\mathcal{A}(x_k,\epsilon_k)$. Hence, $x_k+\gamma(\epsilon_k)s_k^*$ is strictly feasible.

Similarly, we have that with $\gamma = \gamma(\epsilon_k)$ the objective function verifies 
\begin{equation*}
\begin{aligned}   
f_0(x_k+\gamma s_k^*)< f_0(x_k) + 2(M_{\max}+L_{\max}) \gamma ^2 -\epsilon_k \gamma.
\end{aligned}
\end{equation*}
Thus, $f_0(x_k+\gamma s_k^*)<f_0(x_k)-\epsilon^2_k/(8(M_{\max}+L_{\max}))$. 

\section{Proof of Theorem \ref{thm: kkt}}
\label{sec: kkt}
We only consider the case where $\mathcal{A}(x_c,0)$ is not empty. The proof can be easily adapted for $\mathcal{A}(x_c,0)=\emptyset$. 

We show the result through contradiction by assuming that there does not exist $\lambda_c\in \mathbb{R}_{\geq 0}^m$ such that $(x_c,\lambda_c)$ is a KKT pair. Then, one and only one of the following arguments holds:
\begin{itemize}
\item[1)] $\nabla f_0(x_c)$ is not a linear combination of $\nabla f_i(x_c)$, $i\in \mathcal{A}(x_c,0)$,
\item[2)] $\nabla f_0(x_c) = \sum_{i\in\mathcal{A}(x_c,0)} \lambda_i \nabla f_i(x_c)$ and there exists $i^*\in\mathcal{A}(x_c,0)$ such that $\lambda_{i^*}>0$.
\end{itemize}
We show in the following that no matter which argument holds, we can always find $s\in \mathbb{R}^d$ such that 
\begin{equation}
\label{eq: limit_inequalities}
 \langle \nabla f_0(x_c), s \rangle<0,\;\;
\langle \nabla f_i(x_c), s \rangle \leq 0, \;\forall i\in\mathcal{A}(x_c,0).  
\end{equation}

If 1) holds, we let $g_{\parallel}$ be the projection of $\nabla f_0(x_c)$ onto $\textrm{span}\{\nabla f_i(x_c),i\in \mathcal{A}(x_c,0)\}$ and $g_{\perp}:=\nabla f_0(x_c)-g_{\parallel}$. Then $g_\perp\neq 0$, $\langle \nabla f_0(x_c),g_\perp\rangle>0$ and $\langle \nabla f_i(x_c), g_\perp \rangle=0, \;\forall i\in\mathcal{A}(x_c,0)$. Therefore,  $s= -g_\perp$ satisfies \eqref{eq: limit_inequalities}.

If 2) holds, we assume without loss of generality that $\mathcal{A}(x_c,0)\neq \{i^*\}$. Then we let $h_\parallel$ be the projection of $\nabla f_{i^*}(x_c)$ onto $\textrm{span}\{\nabla f_i(x_c),i\in \mathcal{A}(x_c,0)\text{ and }i\neq i^*\}$ and $h_\perp: = \nabla f_{i^*}(x_c) - h_\parallel$. Due to LICQ, $h_\perp\neq 0$. One can verify that  $s = -h_\perp$ also satisfies \eqref{eq: limit_inequalities}.

Then we notice that since the set $\{s:$ \eqref{eq: limit_inequalities} holds$\}$ is non-empty, there exist $\epsilon>0$ and $s_\epsilon$ with $\|s_\epsilon\|_1 = 1$ such that 
\begin{equation}
\label{eq: limit_inequalities_strict}
\langle \nabla f_i(x_c), s_\epsilon\rangle \leq -4\epsilon, \;\forall i\in\mathcal{A}(x_c,0)\cup\{0\}. 
\end{equation}
To see this result, we assume  $\bar{s}\in\mathbb{R}^d$ satisfies  $(\nabla f_0(x_{c}))^\top \bar{s}<0$ and $(\nabla f_i(x_{c}))^\top \bar{s}\leq0$ for any $i\in \mathcal{A}(x_c,0)$. We let $\mathcal{A}(x_c,0)=\{i_1,\ldots,i_l\}$. There exists $y\in \mathbb{R}^d$ such that 
\begin{equation}
\label{eq: steering}
Jy = \begin{bmatrix} -1\\ \vdots \\ -1\end{bmatrix}\text{, where } J=\begin{bmatrix}\nabla f_{i_1}(x_c)^\top \\
\vdots \\ \nabla f_{i_l}(x_c)^\top\end{bmatrix},
\end{equation}
because $J$ is full row rank due to LICQ. Therefore, there exists $\sigma>0$ such that 
\[
\begin{aligned}
\delta_i:=-(\nabla f_i(x_{c}))^\top (\bar{s}+\sigma y)>0,\; \forall i\in \mathcal{A}(x_c,0)\cup \{0\}.
\end{aligned}
\]
Then $s_\epsilon=s^*_\epsilon :=(\bar{s}+\sigma y)/\|\bar{s}+\sigma y\|_1$ and $\epsilon=\epsilon^*:= \frac{1}{4} \min_i{\delta_i}/\|\bar{s}+\sigma y\|_1$ satisfy \eqref{eq: limit_inequalities_strict}.

Due to the continuity of $\nabla f_i(x)$ for $i\in \mathbb{Z}^m_0$, there exists $r>0$ such that any $x\in \mathcal{B}_r(x_c)$ verifies that 
\begin{equation}
\label{eq: within_the_ball}
\langle \nabla f_i(x), s^*_\epsilon\rangle \leq -3\epsilon^*, \;\forall i\in\mathcal{A}(x_c,0)\cup\{0\}. 
\end{equation}
\noindent  
Since $x_c$ is an accmulation point and $\{\epsilon_k\}_{k\geq 1 }$ converges to 0, there exist infinitely many $k$ such that
\begin{equation}
\label{eq: large_enough_k}
\begin{aligned}
&k>K_\mathrm{switch}, \text{ }x_{k}\neq x_{k+1}, \\
&\mathcal{A}(x_c,\epsilon_k)\subset \mathcal{A}(x_c,0),\text{ }x_{k}\in \mathcal{B}_r(x_c).
\end{aligned}
\end{equation}
For any of these $k$s, considering \eqref{eq: within_the_ball} and for any $i$
$$\langle\Delta_i^{\nu^*_k(\epsilon^*)} (x), s^*_\epsilon\rangle \leq |\Delta_i^{\nu^*_k(\epsilon^*)} (x)|\cdot \|s^*_\epsilon\|\leq\epsilon^*,$$
we have 
\begin{equation}
\label{eq: close_to_x_c}
\langle \nabla^{\nu^*_k(\epsilon^*)} f_i(x), s^*_\epsilon\rangle \leq -2\epsilon^*, \;\forall i\in\mathcal{A}(x_c,0)\cup\{0\}. 
\end{equation}

\noindent From Algorithm \ref{alg: SZO-LP}, we see that, for any $k$ satisfying \eqref{eq: large_enough_k}, $x_{k+1}$ is derived through Line 12 and $s^*_k$ through Line 7. Therefore, we can use Proposition \ref{prop: lower_bound_for_epsilon} and \eqref{eq: close_to_x_c} to conclude that $\epsilon_k\geq \epsilon^*/8$ for infinitely many $k$. However, this conclusion contradicts with Point 3 of Theorem \ref{thm: tools_for_convergence}.

\hfill $\blacksquare$}

\end{document}